# THE ASYMPTOTIC DISTRIBUTIONS OF THE LARGEST ENTRIES OF SAMPLE CORRELATION MATRICES

By Tiefeng Jiang

*University of Minnesota*

Let $X_n = (x_{ij})$ be an $n$ by $p$ data matrix, where the $n$ rows form a random sample of size $n$ from a certain $p$-dimensional population distribution. Let $R_n = (\rho_{ij})$ be the $p \times p$ sample correlation matrix of $X_n$; that is, the entry $\rho_{ij}$ is the usual Pearson's correlation coefficient between the $i$th column of $X_n$ and $j$th column of $X_n$. For contemporary data both $n$ and $p$ are large. When the population is a multivariate normal we study the test that $H_0$: the $p$ variates of the population are uncorrelated. A test statistic is chosen as $L_n = \max_{i \neq j} |\rho_{ij}|$. The asymptotic distribution of $L_n$ is derived by using the Chen–Stein Poisson approximation method. Similar results for the non-Gaussian case are also derived.

**1. Introduction.** Let $X_n = (x_{ij})$ be an $n$ by $p$ data matrix, where the $n$ rows are observations from a certain multivariate distribution and each of $p$ columns is an $n$ observation from a variable of the population distribution. Let $\rho_{ij}$ be the Pearson correlation coefficient between the $i$th and $j$th columns of $X_n$. That is

$$(1.1) \qquad \rho_{ij} = \frac{\sum_{k=1}^n (x_{k,i} - \bar{x}_i)(x_{k,j} - \bar{x}_j)}{\sqrt{\sum_{k=1}^n (x_{k,i} - \bar{x}_i)^2} \cdot \sqrt{\sum_{k=1}^n (x_{k,j} - \bar{x}_j)^2}},$$

where $\bar{x}_i = (1/n) \sum_{k=1}^n x_{k,i}$. Then $R_n := (\rho_{ij})$ is a $p$ by $p$ symmetric matrix. It is called the sample correlation matrix generated by $X_n$.

Suppose the population is a multivariate normal distribution with mean vector $\mu$, covariance matrix $\boldsymbol{\Sigma}$ and correlation coefficient matrix $\mathbf{R}$. When the sample size $n$ and the dimension $p$ are large and comparable, Johnstone [14] studied the test with null hypothesis $H_0: \boldsymbol{\Sigma} = \mathbf{I}$ under assumption that $\mu = \mathbf{0}$, where $\mathbf{I}$ is the identity matrix. The null hypothesis is equivalent to that the population distribution is the product of $p$ univariate standard normal distributions. The test statistic is chosen as the maximum eigenvalue of the sample covariance matrix $X_n' X_n$ according to the method principal component analysis (PCA). It is proved that the asymptotic distribution of the maximum eigenvalue is the Tracy–Widom law.

When both $n$ and $p$ are large we consider the test with null hypothesis

$$(1.2) \qquad H_0 : \mathbf{R} = I.$$

Equivalently, the population distribution is a product of univariate normal distribution $N(\mu_i, \sigma_i^2)$'s for some unknown $\mu_i$'s and unknown $\sigma_i$'s. The difference between this test and the one in [14] mentioned above is that all $\mu_i$'s do not have to be identical and all $\sigma_i$'s do not have to be identical, either. Besides, we do not assume that $\mu_i$'s and $\sigma_i$'s are known. Our test seems to be more natural and practical. The maximum eigenvalue $\lambda_{\max}$ of the sample correlation matrix $R_n$ can be taken as the test statistic according to PCA. But the distribution of $\lambda_{\max}$ is not clear so far, although there is evidence that $\lambda_{\max}$ may also follow the Tracy–Widom law asymptotically as shown in [13].

In this paper we do not pursue the maximum eigenvalue $\lambda_{\max}$ as the test statistic because of its complexity. Instead we choose the following intuitive one:

$$L_n = \max_{1 \leq i < j \leq p} |\rho_{ij}|,$$

where $\rho_{ij}$ is as in (1.1). Barbour and Eagleson [6] provided a general idea of dealing with the tail of $L_n$ by using the Poisson approximation method. In this paper we will derive the strong law and limiting distribution of $L_n$ via this method. In fact, we will prove more general results; the observations $x_{ij}$'s do not have to be Gaussian. Our results will be precisely stated next.

Suppose $\{\xi, x_{ij},\ i,j = 1,2,\dots\}$ are i.i.d. random variables. Let $X_n = (x_{ij})_{1 \leq i \leq n, 1 \leq j \leq p}$. Let $x_1, x_2, \dots, x_p$ be the $p$ columns of $X_n$. Then $X_n = (x_1, x_2, \dots, x_p)$. Let $\bar{x}_k$ be the sample average of $x_k$, that is, $\bar{x}_k = (1/n) \sum_{i=1}^n x_{ik}$. We write $x_i - \bar{x}_i$ for $x_i - \bar{x}_i e$, where $e = (1, 1, \dots, 1)^T \in \mathbb{R}^n$. Then, $\rho_{ij}$, the Pearson correlation coefficient in (1.1) between $x_i$ and $x_j$ can be rewritten as

$$(1.3) \qquad \rho_{ij} = \frac{(x_i - \bar{x}_i)^T (x_j - \bar{x}_j)}{\|x_i - \bar{x}_i\| \cdot \|x_j - \bar{x}_j\|}, \qquad 1 \leq i, j \leq p,$$

where $\|\cdot\|$ is the usual Euclidean norm. Obviously, $\rho_{ii} = 1$ for each $i$.

First, we obtain a strong limit theorem as follows.

THEOREM 1.1. *Suppose $E|\xi|^{30-\varepsilon} < \infty$ for any $\varepsilon > 0$. If $n/p \to \gamma \in (0, \infty)$, then*

$$\lim_{n \to \infty} \sqrt{\frac{n}{\log n}} L_n = 2 \qquad a.s.$$

The above strong law of $L_n$ does not depend on $p$ although $X_n$ is an $n$ by $p$ matrix. For the limiting distribution the following holds.

THEOREM 1.2. *Suppose that $E|\xi|^{30+\varepsilon} < \infty$ for some $\varepsilon > 0$. If $n/p \to \gamma$, then*

$$P(nL_n^2 - 4\log n + \log(\log n) \leq y) \to e^{-Ke^{-y/2}}$$

*as $n \to \infty$ for any $y \in \mathbb{R}$, where $K = (\gamma^2 \sqrt{8\pi})^{-1}$.*

The limiting distribution appearing in Theorem 1.2 is called the extreme distribution of type I.

For constants $a_i \in \mathbb{R}^1$ and $b_i \in \mathbb{R}^1$, $i = 1, 2, \ldots, p$, it is easy to see that the matrix $(a_1 x_1 + b_1 e, a_2 x_2 + b_2 e, \ldots, a_p x_p + b_p e)$ and $X_n = (x_1, x_2, \ldots, x_p)$ generate the same correlation matrix $R_n$. Also, if $\xi \sim N(0, 1)$, then $Ee^{t\xi^2} < \infty$ for all $t < 1/2$. We immediately have the following result.

COROLLARY 1.1. *Suppose $\{x_{ij}; i \geq 1, j \geq 1\}$ are independent and $x_{ij} \sim N(\mu_j, \sigma_j^2)$ for some $\mu_j$ and $\sigma_j \neq 0$ for all $i$ and $j$. Let the sample correlation matrix $R_n$ be obtained from $X_n := (x_{ij}; 1 \leq i \leq n, 1 \leq j \leq p)$. Then the conclusions of Theorems 1.1 and 1.2 also hold.*

The above corollary gives the distribution of the test statistic $L_n$ under the null hypothesis in (1.2).

Theorem 1.3 below is used in the proof of Lemmas 3.1 and 3.2. These two lemmas are key to prove Theorems 1.1 and 1.2. It is a nonasymptotic inequality on the moderate deviation of partial sums of independent random variables. Though sums of independent random variables are well understood, we did not notice a similar result in the literature, for example, [19] and [20]. The usual moderate results such as those in [15] and Theorem 3.7.1 on page 109 from [8] are not applicable in our case. The reason is that we do not have identical distribution assumption. Second, asymptotic bounds do not work in our proof because our case involves an uniform bound of infinitely many such probabilities. This is evident from Lemma 2.1 in Section 3. There is a similar situation in the large deviation case. The Chernoff bound (see, e.g., (c) of Remarks on page 27 from [8]) is a nonasymptotic bound of sums of i.i.d. random variables. But the classical Cramér-type large deviation is a limiting result. The Chernoff bound is used in the proof of theorems in [10] and [11] for the same reason of proving our Theorems 1.1 and 1.2 via the following Theorem 1.3.

THEOREM 1.3. *Let $\{\eta_i, 1 \leq i \leq n\}$ be independent random variables with mean zero. Assume $\max_{1 \leq i \leq n} E|\eta_i|^\beta < \infty$ for some $\beta > 2$. Then for any $\rho > 0$* *and* $t > 0$,

$$P\left(\frac{1}{\sqrt{n}}\sum_{i=1}^{n}\eta_i \geq t\right) \leq \frac{M_\beta}{n^{\rho\beta-1}} + K_n e^{-t_\rho^2/(2M_2)}, \tag{1.4}$$

*where*

$$M_s := \frac{1}{n}\sum_{i=1}^{n}E|\eta_i|^s \text{ for } s \in (0,\beta] \quad \text{and} \quad K_n := \exp\left\{\frac{t^3 n^{\rho-(1/2)}}{3M_2^2}e^{2tn^{\rho-(1/2)}/M_2}\right\}$$

*and* $t_\rho = t - M_\beta n^{\rho(1-\beta)+(1/2)}$.

In our applications, $K_n \sim 1$, $t_\rho \sim t$ and $M_2 = 1$. Also, $M_\beta/n^{\rho\beta-1}$ in (1.4) is smaller than the term next to it. So the probability is roughly bounded by $e^{-t^2/2}$.

The main tool used in proving Theorems 1.1 and 1.2 are the Chen–Stein Poisson approximation method and probabilities of moderate deviations by Amosova [1] and Rubin and Sethuraman [21]. They are listed in the Appendix.

In traditional random matrix theories, eigenvalues are the primary concern. See, for example, [18] and [5]. This paper together with [12], in which the maxima of entries of certain Haar-distributed matrices were studied for an imaging analysis problem, suggests that the study of entries of matrices are also important.

Now we state the outline of this paper. A couple of lemmas are given in Section 2 for the preparation of the proofs of main results. We prove all main results in Section 3. In the last section some known results used in the proofs of our theorems are listed.

**2. Auxiliary lemmas.** Three lemmas are needed before we go to the proof of main results. The proof of the following relies on Theorem 1.3, which will be proved at the end. There is no circular reasoning.

LEMMA 2.1. *Let* $\{\xi_k, \eta_k, \eta'_k, k = 1, 2, \ldots, n\}$ *be i.i.d. random variables with mean* 0 *and variance* 1. *Let* $\{u_n; n \geq 1\}$ *be a sequence of positive numbers such that* $u_n/\sqrt{n \log n} \to a \in (0, \infty)$. *If* $E|\xi_1|^q < \infty$ *for some* $q > (a^2+1)(a^2+2)$, *then*

$$P\left(\left|\sum_{k=1}^{n}\xi_k\eta_k\right| \geq u_n, \left|\sum_{k=1}^{n}\xi_k\eta'_k\right| \geq u_n\right) = O(n^{b-a^2}) \tag{2.1}$$

*as* $n \to \infty$ *for any* $b > 0$.

PROOF. The two events in (2.1) are conditionally independent given $\xi_k$'s. Denote by $P^1$ and $E^1$ such conditional probability and expectation, respectively. Then the probability in (2.1) is

$$(2.2) \qquad E\left[P^1\left(\left|\sum_{k=1}^n \xi_k \eta_k\right| \geq u_n\right)^2\right].$$

Set

$$A_n(s) = \left\{\frac{1}{n}\left|\sum_{k=1}^n (|\xi_k|^s - E|\xi_k|^s)\right| \leq \delta\right\}$$

for $s \geq 2$ and $\delta \in (0, 1/2)$. Choose $\beta \in (a^2+2, q/(a^2+1))$ and $r = a^2+1$. Let $\zeta_k = |\xi_k|^\beta - E|\xi_k|^\beta$ for $k = 1, 2, \ldots, n$. Then $E|\zeta_1|^r < \infty$. By the Chebyshev inequality and Lemma A.1,

$$(2.3) \quad P(A_n(\beta)^c) = P\left(\left|\sum_{k=1}^n \zeta_k\right| > n\delta\right) \leq (n\delta)^{-r} E\left|\sum_{k=1}^n \zeta_k\right|^r = O(n^{-f(r)})$$

as $n \to \infty$, where $f(r) = r/2$ if $r \geq 2$, and $f(r) = r-1$ if $1 < r \leq 2$. Let $\{\zeta_k'; 1 \leq k \leq n\}$ be an independent copy of $\{\zeta_k; 1 \leq k \leq n\}$. Then since, from (2.3), $P(|\sum_{k=1}^n \zeta_k| \leq n\delta/2) \geq 1/2$ for sufficiently large $n$, it follows that

$$(2.4) \quad P\left(\left|\sum_{k=1}^n \zeta_k\right| > n\delta\right) \leq 2P\left(\left|\sum_{k=1}^n (\zeta_k - \zeta_k')\right| > n\delta/2\right) = O(n^{-f(r)})$$

by repeating (2.3). Given an integer $j \geq 1$, let $v = n\delta/4j$. Then by Lemma A.2, there are positive constants $C_j$ and $D_j$ such that

$$P\left(\left|\sum_{k=1}^n (\zeta_k - \zeta_k')\right| > n\delta/2\right)$$

$$= P\left(\left|\sum_{k=1}^n (\zeta_k - \zeta_k')\right| > 2jv\right)$$

$$\leq C_j P\left(\max_{1 \leq k \leq n} |\zeta_k - \zeta_k'| > v\right) + D_j P\left(\left|\sum_{k=1}^n (\zeta_k - \zeta_k')\right| > v\right)^j.$$

Since $E|\zeta_1|^r < \infty$, $P(\max_{1 \leq k \leq n} |\zeta_k - \zeta_k'| > v) \leq nP(|\zeta_1 - \zeta_1'| > v) = O(n^{1-r})$. By the same argument as the equality in (2.4), we obtain

$$\left(P\left(\left|\sum_{k=1}^n (\zeta_k - \zeta_k')\right| > v\right)\right)^j = O(n^{-jf(r)}).$$

Take $j = [(r-1)/f(r)] + 1$. It follows that

$$(2.5) \qquad P\left(\left|\sum_{k=1}^n (\zeta_k - \zeta_k')\right| > n\delta/2\right) = O(n^{1-r})$$

as $n \to \infty$. Combining (2.3) and (2.4) with (2.5), we obtain that
$$P(A_n(\beta)^c) = O(n^{1-r})$$
as $n \to \infty$. By the same arguments the above still holds if $\beta$ is replaced by 2. Consequently,

$$
\begin{aligned}
(2.6) \quad & E\left[P^1\left(\left|\sum_{k=1}^n \xi_k \eta_k\right| \geq u_n\right)^2\right] \\
& \leq E\left[P^1\left(\left|\sum_{k=1}^n \xi_k \eta_k\right| \geq u_n\right)^2 I_{A_n(2) \cap A_n(\beta)}\right] + O(n^{1-r}).
\end{aligned}
$$

Now we apply Theorem 1.3 to the last probability in (2.6). Note $E^1(\xi_k \eta_k) = 0$ and $E^1|\xi_k \eta_k|^s = |\xi_k|^s (E|\xi_1|^s)$ for any $s > 0$. In particular, $E^1(\xi_k \eta_k)^2 = \xi_k^2$. Thus
$$M_s = \frac{1}{n}\left(\sum_{k=1}^n |\xi_k|^s\right) \cdot E|\xi_1|^s.$$

So $1 - \delta < M_2 \leq 1 + \delta$ on $A_n(2)$. Since $\delta \in (0, 1/2)$, $1/2 \leq M_2 \leq 2$ on $A_n(2)$. Moreover, $M_\beta \leq (1 + E|\xi_1|^\beta)E|\xi_1|^\beta < \infty$ on $A_n(\beta)$. Choose $t = u_n/\sqrt{n}$ and $\rho \in ((a^2 + 2)/(2\beta), 1/2)$. It is easy to verify that $d := -\rho(1 - \beta) - (1/2) > 0$,
$$|t - t_\rho| \leq \frac{(1 + E|\xi_1|^\beta)E|\xi_1|^\beta}{n^d} \quad \text{and} \quad K_n \leq \exp(2t^3 n^{\rho - 1/2} e^{4tn^{\rho-1/2}})$$
on $A_n(2) \cap A_n(\beta)$ for each $n \geq 1$. Then there is a constant $C > 0$, such that the probability in (2.6) under the restriction $A_n(2) \cap A_n(\beta)$ is less than
$$C(n^{-a^2/2} + n^{-a^2/(2(1+\delta))}) = (n^{-a^2/(2+2\delta)})$$
for $n$ sufficiently large, where the fact $\rho\beta > 1 + (a^2/2)$ is used. Note that $O(n^{1-r}) = O(n^{-a^2})$ since $r = 1 + a^2$. So the left-hand side of (2.6), hence, the probability in (2.1) is $O(n^{-a^2/(1+\delta)})$ by (2.2). The desired conclusion then follows by choosing $\delta$ small enough. $\square$

For any square matrix $A = (a_{i,j})$, define $\|A\| = \max_{1 \leq i \neq j \leq n} |a_{i,j}|$; that is, the maximum of the absolute values of the off-diagonal entries of $A$.

LEMMA 2.2. *Recall $x_i$ in (1.3). Let $h_i = \|x_i - \bar{x}_i\|/\sqrt{n}$ for each $i$. Then*
$$\|nR_n - X_n^T X_n\| \leq (b_1^2 + 2b_1)W_n b_3^{-2} + nb_3^{-2} b_4^2,$$
*where*
$$b_1 = \max_{1 \leq i \leq n} |h_i - 1|, \qquad W_n = \max_{1 \leq i < j \leq n} |x_i^T x_j|,$$
$$b_3 = \min_{1 \leq i \leq n} h_i, \qquad b_4 = \max_{1 \leq i \leq n} |\bar{x}_i|.$$

PROOF. As in (1.3), the $(i,j)$-entry of $R_n$ is

$$\rho_{ij} = \frac{(x_i - \bar{x}_i)^T(x_j - \bar{x}_j)}{\|x_i - \bar{x}_i\| \cdot \|x_j - \bar{x}_j\|} = \frac{x_i^T x_j - n\bar{x}_i \bar{x}_j}{n h_i h_j}.$$

The $(i,j)$-entry of $X_n^T X_n$ is $x_i^T x_j$. So

$$|n\rho_{ij} - x_i^T x_j| \leq |(h_i h_j)^{-1} - 1| \cdot |x_i^T x_j| + n\frac{|\bar{x}_i|}{h_i} \cdot \frac{|\bar{x}_j|}{h_j}.$$

Taking maximum for both sides, we obtain

$$\|nR_n - X_n^T X_n\|$$
$$\leq \max_{1 \leq i < j \leq n} |(h_i h_j)^{-1} - 1| \cdot \max_{1 \leq i < j \leq n} |x_i^T x_j| + n\left(\max_{1 \leq i \leq n} \frac{|\bar{x}_i|}{h_i}\right)^2.$$

Write $1 - (h_i h_j)^{-1} = (h_i h_j)^{-1}((h_i - 1)(h_j - 1) + (h_i - 1) + (h_j - 1))$. Then the desired inequality follows. □

Next we estimate $b_i$'s.

LEMMA 2.3. *Suppose that $\{\xi, x_{ij}, i, j = 1, 2, \ldots\}$ are i.i.d. random variables with $E\xi = 0$ and $\text{Var}(\xi) = 1$. Suppose also $n/p \to \gamma \in (0, \infty)$. If $E|\xi|^{4/(1-\alpha)} < \infty$ for some $\alpha \in (0, 1/2)$, then*

$$n^\alpha b_1 \to 0 \quad a.s., \quad b_3 \to 1 \quad a.s. \quad \text{and} \quad n^\alpha b_4 \to 0 \quad a.s.$$

*as $n \to \infty$.*

PROOF. The second limit follows from the first one. Easily, $\|x_i - \bar{x}_i\|^2 = x_i^T x_i - n|\bar{x}_i|^2$. Using the fact that $|x - 1| \leq |x^2 - 1|$ for any $x > 0$, we have that

$$(2.7) \qquad n^\alpha b_1 \leq \max_{1 \leq i \leq n}\left|\frac{x_i^T x_i - n}{n^{1-\alpha}}\right| + \left(n^{\alpha/2} \max_{1 \leq i \leq n} |\bar{x}_i|\right)^2.$$

Note $x_i^T x_i = \sum_{k=1}^n x_{ki}^2$. By Lemma A.5 the first and the second maxima above go to zero when $E|\xi|^{4/(1-\alpha)} < \infty$. So the first limit is proved. Under the condition that $E|\xi|^{2/(1-\alpha)} < \infty$, the limit that $n^\alpha b_4 \to 0$ a.s. is proved by noting the relationship between $n^\alpha b_4$ and the right most term in (2.7). □

The analysis of $W_n$ is given in the next section.

**3. Proof of main results.** Recall the definition of $\rho_{ij}$ in (1.1) and (1.3). To prove Theorems 1.1 and 1.2, we assume throughout this section, without loss of generality, that

$$\{\xi, x_{ij}; i, j = 1, 2, \ldots\} \text{ are i.i.d. with } E\xi = 0 \quad \text{and} \quad \text{Var}(\xi) = 1.$$

The proofs of Theorems 1.1 and 1.2 rely on an analysis of the covariance matrix $X_n^T X_n$. The $(i, j)$-entry of $X_n^T X_n$ is $\sum_{k=1}^n x_{ki} x_{kj}$. Recall

$$(3.1) \qquad W_n = \max_{1 \leq i < j \leq n} \left| \sum_{k=1}^n x_{ki} x_{kj} \right|$$

as in Lemma 2.2. The first step in proving our main theorems is approximating $R_n$ by $X_n^T X_n$ as shown in Lemmas 2.2 and 2.3. The second step is deriving the corresponding results for $X_n^T X_n$. We actually will prove the following two lemmas.

LEMMA 3.1. *Suppose that $E|\xi|^{30-\varepsilon} < \infty$ for any $\varepsilon > 0$. If $n/p \to \gamma \in (0, \infty)$, then:*

(i) $\limsup_{n \to \infty} \dfrac{W_n}{\sqrt{n \log n}} \leq 2 \quad a.s.$

(ii) $\liminf_{n \to \infty} \dfrac{W_n}{\sqrt{n \log n}} \geq 2 \quad a.s.$

Lemma 3.1 actually says that $W_n/\sqrt{n \log n} \to 2$ a.s. as $n \to \infty$. The reason we did not combine (i) and (ii) as a single limit is that the proof of the combined one is relatively long. We will prove the two parts separately.

LEMMA 3.2. *Suppose that $E|\xi|^{30+\varepsilon} < \infty$ for some $\varepsilon > 0$. If $n/p \to \gamma \in (0, \infty)$, then*

$$P\left(\frac{W_n^2 - \alpha_n}{n} \leq y\right) \to e^{-Ke^{-y/2}}$$

*as $n \to \infty$ for any $y \in \mathbb{R}$, where $\alpha_n = 4n \log n - n \log(\log n)$ and $K = (\gamma^2 \sqrt{8\pi})^{-1}$.*

Assuming Lemmas 3.1 and 3.2, we next prove Theorems 1.1 and 1.2. The proof of the former two lemmas are given later.

PROOF OF THEOREMS 1.1 AND 1.2. Choose $\alpha = 1/3$. Under the condition that $E|\xi|^6 < \infty$, we have from the triangle inequality, Lemmas 2.2 and 2.3 that

$$(3.2) \quad |nL_n - W_n| \leq \|nR_n - X_n^T X_n\| \leq 4n^{-1/3} W_n + 2n^{1/3} \qquad a.s.$$

as $n$ is sufficiently large. Applying Lemma 3.1, it follows that $4n^{-1/3}W_n = O(n^{1/6}\log n)$ almost surely. Hence $nL_n - W_n = O(n^{1/3})$ a.s. Theorem 1.1 then follows immediately from Lemma 3.1. Now Theorem 1.1 and Lemma 3.1 imply that $nL_n + W_n = O(\sqrt{n\log n})$. Consequently,

$$nL_n^2 - \frac{W_n^2}{n} = \frac{1}{n}(nL_n - W_n)(nL_n + W_n)$$
$$= O(n^{-1/6}(\log n)^{1/2}) \quad \text{a.s.}$$

Theorem 1.2 then follows from Lemma 3.2. □

Now we turn to prove Lemmas 3.1 and 3.2.

PROOF OF LEMMA 3.1(i). Given $\delta \in (0,1)$, let $w_n = (2+\delta)\sqrt{n\log n}$. Define $y_{ij}^{(l)} := \sum_{k=1}^{l} x_{ki}x_{kj}$, $i,j,l \geq 1$. Then $y_{ij}^{(l)}$ is a sum of $l$ i.i.d. random variables with mean zero and variance one. By Lemma A.3, under the condition that $E|\xi|^d < \infty$ for some $d > 2 + (2+\delta)^2$,

$$(3.3) \qquad \max_{1 \leq i \neq j < \infty} P(|y_{ij}^{(l)}| > w_l) = O\left(\frac{1}{l^{2+\delta}}\right)$$

as $l$ is large, where we also use the fact that

$$(3.4) \qquad 1 - \Phi(x) = \frac{1}{\sqrt{2\pi}} \int_x^\infty e^{-t^2/2}\,dt \sim \frac{1}{\sqrt{2\pi}\,x} e^{-x^2/2}$$

as $x \to +\infty$ (see, e.g., page 49 from [7]). Review the expression of $W_n$ in (3.1). For any integer $m > 4/\delta$,

$$(3.5) \qquad \begin{aligned}\max_{n^m \leq l \leq (n+1)^m} W_l &\leq \max_{1 \leq i \neq j \leq (n+1)^m}\left(\max_{n^m \leq l \leq (n+1)^m} |y_{ij}^{(l)}|\right) \\ &\leq \max_{1 \leq i \neq j \leq (n+1)^m} |y_{ij}^{(n^m)}| + r_n,\end{aligned}$$

where

$$(3.6) \qquad r_n = \max_{1 \leq i \neq j \leq (n+1)^m} \max_{n^m \leq l \leq (n+1)^m} |y_{ij}^{(l)} - y_{ij}^{(n^m)}|.$$

By (3.3),

$$P\left(\max_{1 \leq i \neq j \leq (n+1)^m} |y_{ij}^{(n^m)}| > w_{n^m}\right) \leq (n+1)^{2m} P(|y_{12}^{(n^m)}| > w_{n^m})$$
$$= O(n^{-\delta m}).$$

Since $\sum_n n^{-\delta m} < \infty$, by the Borel–Cantelli lemma,

$$(3.7) \qquad \limsup_{n \to \infty} \frac{\max_{1 \leq i \neq j \leq (n+1)^m} |y_{ij}^{(n^m)}|}{\sqrt{n^m \log(n^m)}} \leq 2+\delta \quad \text{a.s.}$$

Now let us estimate $r_n$ as in (3.6).

Let $\{z_1, z_2, \ldots\}$ be i.i.d. random variables with the same law as $x_{11}x_{12}$ with partial sums $S_0 = 0$ and $S_k = \sum_{i=1}^k z_i$. Clearly, $Ez_1 = 0$ and $Ez_1^2 = 1$. Observe that the distribution of $y_{ij}^{(l)} - y_{ij}^{(n^m)}$ is equal to that of $S_{l-n^m}$ for all $l \geq n^m$. Thus,

$$P(r_n \geq \delta\sqrt{n^m \log(n^m)})$$

$$(3.8) \qquad \leq (n+1)^{2m} P\left(\max_{1 \leq k \leq (n+1)^m - n^m} |S_k| \geq \delta\sqrt{n^m \log(n^m)}\right)$$

$$\leq 2(n+1)^{2m} P(|S_{(n+1)^m - n^m}| \geq (\delta/2)\sqrt{n^m \log(n^m)})$$

as $n$ is sufficiently large, where Ottaviani's inequality (see Exercise 16 on page 74 in [7]) is used in the last inequality. Set $k_n = (n+1)^m - n^m$. Note that, for fixed $m$ and $\delta$, $(\delta/2)\sqrt{n^m \log(n^m)} \geq (2+\delta)\sqrt{k_n \log k_n}$ as $n$ is sufficiently large. By (3.4) and Lemma A.3, the last probability in (3.8) is equal to $O(\exp(-(m-1)(2+\delta)^2(\log n)/2))$ provided $E|\xi|^d < \infty$ for some $d > 2 + (2+\delta)^2$. Therefore,

$$P(r_n \geq \delta\sqrt{n^m \log(n^m)}) = O(n^{-u}),$$

where $u = (m-1)(2+\delta)^2/2 - 2m > 1$, since $m$ is chosen such that $m > 4/\delta$. By the Borel–Cantelli lemma again,

$$(3.9) \qquad \limsup_{n \to \infty} \frac{r_n}{\sqrt{n^m \log(n^m)}} \leq \delta \qquad \text{a.s.}$$

By (3.5), (3.7) and (3.9), we obtain that

$$\limsup_{n \to \infty} \frac{\max_{n^m \leq l \leq (n+1)^m} W_l}{\sqrt{n^m \log(n^m)}} \leq 2 + 2\delta \qquad \text{a.s.}$$

for any sufficiently small $\delta > 0$. This implies inequality (i) in Lemma 3.1. □

PROOF OF LEMMA 3.1(ii). We continue to use the notations in the proof of (i) of Lemma 3.1. For any $\delta \in (0,1)$, define $v_n = (2-\delta)\sqrt{n \log n}$. We first claim that

$$(3.10) \qquad P(W_n \leq v_n) = O\left(\frac{1}{n^{\delta'}}\right)$$

as $n \to \infty$ for some positive constant $\delta'$ depending on $\delta$ and the distribution of $\xi$ only. If this is true, take an integer $m$ such that $m > 1/\delta'$. Then $P(W_{n^m} \leq v_{n^m}) = O(1/n^{\delta'm})$. Since $\sum_n n^{-\delta'm} < \infty$, by the Borel–Cantelli lemma, we have that

$$(3.11) \qquad \liminf_{n \to \infty} \frac{W_{n^m}}{\sqrt{n^m \log(n^m)}} \geq 2 - \delta \qquad \text{a.s.}$$

for any $\delta \in (0,1)$. Recalling the definition of $r_n$ in (3.6), we have that

$$\inf_{n^m \leq k \leq (n+1)^m} W_k \geq W_{n^m} - r_n.$$

By (3.9) and (3.11), we have that

$$\liminf_{n \to \infty} \frac{\inf_{n^m \leq k \leq (n+1)^m} W_k}{\sqrt{n^m \log(n^m)}} \geq 2 - 2\delta \quad \text{a.s.}$$

for any $\delta$ small enough. This implies (ii) of Lemma 3.1.

Now we turn to prove claim (3.10) by Lemma A.4.

Take $I = \{(i,j); \ 1 \leq i < j \leq p\}$. For $\alpha = (i,j) \in I$, set $B_\alpha = \{(k,l) \in I;$ one of $k$ and $l = i$ or $j$, but $(k,l) \neq \alpha\}$, $\eta_\alpha = |y_{ij}^{(n)}|$, $t = v_n$ and $A_\alpha = A_{ij} = \{|y_{ij}^{(n)}| > v_n\}$. By Lemma A.4,

(3.12) $$P(W_n \leq v_n) \leq e^{-\lambda_n} + b_{1,n} + b_{2,n}.$$

Evidently

(3.13)
$$\lambda_n = \frac{p(p-1)}{2} P(A_{12}),$$
$$b_{1,n} \leq 2p^3 P(A_{12})^2 \quad \text{and} \quad b_{2,n} \leq 2p^3 P(A_{12} A_{13}).$$

Remember that $y_{12}^{(n)}$ is a sum of i.i.d. random variables with mean 0 and variance 1. Recall (3.4). By Lemma A.3,

(3.14) $$P(A_{12}) \sim \frac{1}{(2-\delta)\sqrt{2\pi \log n}} \cdot \frac{1}{n^{(2-\delta)^2/2}}$$

as $n \to \infty$ provided $E|\eta|^6 < \infty$. Note that $P(A_{12} A_{13}) = P(|y_{12}^{(n)}| \geq v_n, |y_{13}^{(n)}| \geq v_n)$ and $v_n/\sqrt{n \log n} \to 2 - \delta$. By Lemma 2.1, $P(A_{12} A_{13}) = O(n^{b-(2-\delta)^2})$ for any $b > 0$ provided $E|\xi|^q < \infty$ for some $q > ((2-\delta)^2 + 1)((2-\delta)^2 + 2) < 30$. Choosing both $b$ and $\delta$ small enough, we obtain

(3.15) $$e^{-\lambda_n} \leq e^{-n^\delta}, \qquad b_{1,n} \leq \frac{1}{\sqrt{n}} \quad \text{and} \quad b_{2,n} \leq \frac{1}{\sqrt{n}}$$

for sufficiently large $n$. Then (3.10) follows from (3.12) and (3.15). $\square$

PROOF OF LEMMA 3.2. We need to show that

(3.16) $$P\left(\max_{1 \leq i < j \leq p} |y_{ij}| \leq \sqrt{\alpha_n + ny}\right) \to e^{-Ke^{-y/2}},$$

where $y_{ij} = \sum_{k=1}^n x_{ki} x_{kj}$. Now we apply Lemma A.4 to prove (3.16). Take $I = \{(i,j); 1 \leq i < j \leq p\}$. For $\alpha = (i,j) \in I$, set $X_\alpha = |y_{ij}|$ and $B_\alpha = \{(k,l) \in I;$ one of $k$ and $l =$ $i$ or $j$, but $(k,l) \neq \alpha\}$. Choose $t = \sqrt{\alpha_n + ny}$. We first calculate $\lambda = \lambda_n$ in the theorem. Since $\{y_{ij}; (i,j) \in I\}$ are identically distributed,

$$
\begin{aligned}
\lambda_n &= \sum_{1 \leq i < j \leq p} P(|y_{ij}| > \sqrt{\alpha_n + ny}) \\
&= \frac{p^2 - p}{2} P\left(\frac{|y_{12}|}{\sqrt{n}} > \sqrt{\frac{\alpha_n}{n} + y}\right).
\end{aligned}
\tag{3.17}
$$

Observe that $y_{12}$ is a sum of i.i.d. random variables with mean 0 and variance 1. Since $\sqrt{(\alpha_n/n) + y} \sim 2\sqrt{\log n}$ as $n \to \infty$, it follows from Lemma A.3 that

$$
\begin{aligned}
P&\left(\frac{|y_{12}|}{\sqrt{n}} > \sqrt{(\alpha_n/n) + y}\right) \\
&= P\left(\frac{y_{12}}{\sqrt{n}} > \sqrt{(\alpha_n/n) + y}\right) + P\left(\frac{-y_{12}}{\sqrt{n}} > \sqrt{(\alpha_n/n) + y}\right) \\
&\sim \frac{e^{-y/2}}{\sqrt{2\pi}} n^{-2}
\end{aligned}
\tag{3.18}
$$

provided $E|\xi|^q < \infty$ for some $q > c^2 + 2 = 6$. Thus

$$
\lambda_n \to \frac{e^{-y/2}}{\gamma^2 \sqrt{8\pi}}.
\tag{3.19}
$$

Obviously, $X_\alpha$ is independent of $\{X_\beta; \beta \in I \setminus B_\alpha\}$ for any $\alpha = (i,j) \in I$. To complete (3.16), by Lemma A.4, we have to verify that $b_1 \to 0$ and $b_2 \to 0$ as $n \to \infty$. It is easy to check that the size of $B_\alpha$ is less than $2p$. Thus

$$
b_1 \leq \frac{1}{2}(p^2 - p) \cdot (2p) \cdot P\left(\frac{|y_{12}|}{\sqrt{n}} > \sqrt{\frac{\alpha_n}{n} + y}\right)^2 = O\left(\frac{1}{n}\right)
$$

by (3.18). Also, by symmetry,

$$
b_2 \leq p(p^2 - p) P(|y_{12}| > \sqrt{\alpha_n + ny}, \ |y_{13}| > \sqrt{\alpha_n + ny}).
\tag{3.20}
$$

Here $\sqrt{\alpha_n + ny}/\sqrt{n \log n} \to 2$. By Lemma 2.1, the above probability is $O(n^{b-4})$ for any $b > 0$, provided $E|\xi|^q < \infty$ for some $q > (2^2 + 1)(2^2 + 2) = 30$. Now choose $b < 1$, then $b_2 \to 0$. By Lemma A.4, (3.16) is concluded. $\square$

PROOF OF THEOREM 1.3. Define

$$
\tilde{\eta}_i = \eta_i I(|\eta_i| \leq n^\rho) - E\eta_i I(|\eta_i| \leq n^\rho)
$$

for $\rho > 0$. It is easy to see that

$$
P\left(\frac{1}{\sqrt{n}} \sum_{i=1}^n \eta_i \geq t\right) \leq P\left(\left(\frac{1}{\sqrt{n}} \sum_{i=1}^n \tilde{\eta}_i\right) + r_n \geq t\right) + \sum_{i=1}^n P(|\eta_i| \geq n^\rho),
$$

where $r_n = (1/\sqrt{n}) \sum_{i=1}^n E\eta_i I(|\eta_i| > n^\rho)$ because $E\eta_i = 0$. Clearly, by Markov's inequality,

$$\sum_{i=1}^n P(|\eta_i| \geq n^\rho) \leq \frac{M_\beta}{n^{\rho\beta-1}} \quad \text{and} \quad |r_n| \leq \frac{M_\beta}{n^{\rho(\beta-1)-(1/2)}}.$$

Thus, to prove the lemma, it suffices to show that

$$(3.21) \qquad P\left(\frac{1}{\sqrt{n}} \sum_{i=1}^n \tilde{\eta}_i \geq t_\rho\right) \leq K_n e^{-t_\rho^2/(2M_2)}.$$

By the Chebyshev inequality and independence, we obtain that

$$(3.22) \qquad P\left(\frac{1}{\sqrt{n}} \sum_{i=1}^n \tilde{\eta}_i \geq t_\rho\right) \leq e^{-\theta t_\rho} \prod_{i=1}^n E e^{\theta \tilde{\eta}_i/\sqrt{n}}$$

for any $\theta > 0$. Since $e^x \leq 1 + x + (x^2/2) + (|x|^3/6)e^{|x|}$ for any $x \in \mathbb{R}$,

$$(3.23) \qquad E e^{\theta \tilde{\eta}_i/\sqrt{n}} \leq 1 + \frac{\theta^2}{2n} E\tilde{\eta}_i^2 + \frac{\theta^3}{6n^{3/2}} E(|\tilde{\eta}_i|^3 \exp(\theta|\tilde{\eta}_i|/\sqrt{n})).$$

Obviously, $E\tilde{\eta}_i^2 \leq E\eta_i^2$ and $|\tilde{\eta}_i| \leq 2n^\rho$. It follows that

$$E(|\tilde{\eta}_i|^3 \exp(\theta|\tilde{\eta}_i|/\sqrt{n})) \leq 2n^\rho e^{2\theta n^{\rho-(1/2)}} E|\eta_i|^2.$$

Since $1 + x \leq e^x$ for any $x \in \mathbb{R}$, by (3.23), we have that

$$E e^{\theta \tilde{\eta}_i/\sqrt{n}} \leq \exp\left\{\frac{\theta^2}{2n} E(\eta_i^2) + \frac{\theta^3 E|\eta_i|^2}{3n^{(3/2)-\rho}} e^{2\theta n^{\rho-(1/2)}}\right\}.$$

Substituting this back to (3.22), we obtain

$$P\left(\frac{1}{\sqrt{n}} \sum_{i=1}^n \tilde{\eta}_i \geq t_\rho\right) \leq \exp\left(-\theta t_\rho + \frac{M_2}{2}\theta^2 + \frac{\theta^3 M_2}{3n^{(1/2)-\rho}} e^{2\theta n^{\rho-(1/2)}}\right)$$

for any $\theta > 0$. Choosing $\theta = t_\rho/M_2$, it follows that

$$P\left(\frac{1}{\sqrt{n}} \sum_{i=1}^n \tilde{\eta}_i \geq t_\rho\right) \leq e^{-t_\rho^2/(2M_2)} \exp\left\{\frac{t_\rho^3 n^{\rho-(1/2)}}{3M_2^2} e^{2t_\rho n^{\rho-(1/2)}/M_2}\right\}$$

$$\leq K_n e^{-t_\rho^2/(2M_2)}$$

since $t_\rho < t$ and $t > 0$, where $K_n$ is as in the statement of Theorem 1.3. Then (3.21) follows. The proof is complete. $\square$

## APPENDIX

For the proofs of the main theorems we quote some results from literature in this section.

The following is a corollary of the Marcinkiewicz–Zygmund inequality, see, for example, Corollary 2 on page 368 for $p \geq 2$ and Theorem 2 on page 367 for $p \in [1,2)$ from [7].

LEMMA A.1. *If $\{\eta_n, n \geq 1\}$ are i.i.d. random variables with $E\eta_1 = 0$, $E|\eta_1|^p < \infty$, $p \geq 1$, and $S_n = \sum_{i=1}^n \eta_i$. Then,*

$$E|S_n|^p = \begin{cases} O(n^{p/2}), & \text{if } p \geq 2, \\ O(n), & \text{if } 1 \leq p < 2. \end{cases}$$

The following is Lemma 2.2 from [17], which is a useful version of the maximal inequality of Hoffmann-Jøgensen, see [9] or Proposition 6.7 from [16].

LEMMA A.2. *Let $\{\eta_k, 1 \leq k \leq n\}$ be independent symmetric random variables and $S_n = \sum_{k=1}^n \eta_k$. Then, for each integer $j \geq 1$, there exist positive numbers $C_j$ and $D_j$ depending only on $j$ such that for all $t > 0$,*

$$P(|S_n| \geq 2jt) \leq C_j P\left(\max_{1 \leq k \leq n} |\eta_j| \geq t\right) + D_j (P(|S_n| \geq t))^j.$$

The following lemma is from [1]. It is a refinement of Theorem 2 from [21]. See also page 254 from [20].

LEMMA A.3. *Let $\{\eta_1, \eta_2, \ldots\}$ be a sequence of independent random variables with $E\eta_i = 0$ and $E\eta_i^2 = \sigma_i^2 > 0$ for all $i$. Define $S_n = \sum_{i=1}^n \eta_i$ and $B_n = \sum_{i=1}^n \sigma_i^2$. Suppose*

$$\liminf_n \left\{\frac{B_n}{n}\right\} > 0 \quad \text{and} \quad \limsup_n \left\{\frac{1}{n} \sum_{i=1}^n E|\eta_i|^q\right\} < \infty$$

*for some $q > 2$. Let $\Phi(x) = \int_{-\infty}^x (2\pi)^{-1/2} e^{-t^2/2} dt$. Then*

$$\frac{P(S_n/\sqrt{B_n} > x)}{1 - \Phi(x)} \to 1$$

*uniformly on $[0, c\sqrt{\log n}\,]$ for any $c \in (0, \sqrt{q-2}\,)$ as $n \to \infty$.*

The following Poisson approximation result is essentially a special case of Theorem 1 in [3], which is again a special case of the general Chen–Stein Poisson approximation method. One application of the following lemma is studying behaviors of maxima of random variables. See, for example, [10] and [11].

LEMMA A.4. *Let $I$ be an index set and $\{B_\alpha, \alpha \in I\}$ be a set of subsets of $I$, that is, $B_\alpha \subset I$ for each $\alpha \in I$. Let also $\{\eta_\alpha, \alpha \in I\}$ be random variables. For a given $t \in \mathbb{R}$, set $\lambda = \sum_{\alpha \in I} P(\eta_\alpha > t)$. Then*

$$\left| P\left(\max_{\alpha \in I} \eta_\alpha \leq t\right) - e^{-\lambda} \right| \leq (1 \wedge \lambda^{-1})(b_1 + b_2 + b_3),$$

*where*

$$b_1 = \sum_{\alpha \in I} \sum_{\beta \in B_\alpha} P(\eta_\alpha > t) P(\eta_\beta > t),$$

$$b_2 = \sum_{\alpha \in I} \sum_{\alpha \neq \beta \in B_\alpha} P(\eta_\alpha > t, \eta_\beta > t),$$

$$b_3 = \sum_{\alpha \in I} E|P(\eta_\alpha > t | \sigma(\eta_\beta, \beta \notin B_\alpha)) - P(\eta_\alpha > t)|,$$

*and $\sigma(\eta_\beta, \beta \notin B_\alpha)$ is the $\sigma$-algebra generated by $\{\eta_\beta, \beta \notin B_\alpha\}$. In particular, if $\eta_\alpha$ is independent of $\{\eta_\beta, \beta \notin B_\alpha\}$ for each $\alpha$, then $b_3 = 0$.*

The following is Lemma 2 from [4].

LEMMA A.5. *Let $\{v, v_{ij}, i, j = 1, 2, \ldots\}$ be a double array of i.i.d. random variables and let $\alpha > 1/2$, $\beta > 0$ and $M > 0$ be constants. Then as $n \to \infty$,*

$$\max_{1 \leq j \leq Mn^\beta} \left| n^{-\alpha} \sum_{i=1}^n (v_{ij} - Ev) \right| \to 0 \qquad a.s.,$$

*if and only if $E|v|^{(\beta+1)/\alpha} < \infty$.*

**Acknowledgments.** The author thanks Iain Johnstone very much for suggesting the study of the largest entry of the covariance matrix in the beginning of Section 3. The author also thanks an anonymous referee for reading the manuscript carefully and pointing out a relevant paper.

SCHOOL OF STATISTICS
UNIVERSITY OF MINNESOTA
313 FORD HALL
224 CHURCH STREET S.E.
MINNEAPOLIS, MINNESOTA 55455
USA
E-MAIL: tjiang@stat.umn.edu